\documentclass[a4paper,12pt]{amsart}
\usepackage[latin1]{inputenc}
\usepackage[T1]{fontenc}
\usepackage{amsfonts,amssymb}
\usepackage{amsmath}
\usepackage{pstricks}
\usepackage{amsthm}
\usepackage[all]{xy}
\usepackage[pdftex]{graphicx}
\usepackage{graphicx}   
\newcommand{\n}[0]{{\textbf N}}

\renewcommand{\c}[0]{{\textbf C}}

\renewcommand{\k}[0]{{\textbf K}}

\newtheorem{theo}{Theorem}[section]
\newtheorem{defi}{Definition}[section]
\newtheorem{prop}[theo]{Proposition}

\newtheorem{cor}[theo]{Corollary}
\newtheorem{lem}[theo]{Lemma}

\title[Partition identities and Groebner basis]{Partition identities and application to infinite dimensional Groebner basis and viceversa}

\author{Pooneh Afsharijoo, Hussein Mourtada}

\begin {document}

\maketitle  
\begin{abstract} In the first part of this article, we consider a  Groebner basis of the differential ideal $[x_1^2]$ with respect to "the" weighted 
lexicographical monomial order and show that its computation is related with an identity involving the partitions that
appear in the first Rogers-Ramanujan identity. We then prove that a Grobener basis of this ideal is not differentially finite in contrary 
with the case of "the" weighted reverse lexicographical order. In the second part, we give a simple and direct proof of a theorem of Nguyen Duc Tam about the Groebner basis of the differential ideal $[x_1y_1];$ we then obtain identities involving partitions with $2$ colors.

\end{abstract}

\footnote{{\textbf{2010 Mathematics Subject Classification.} 
05A17,12H05,13D40,13P10.}}

\section*{Intoduction}

An integer partition of a positive integer number $n$ is a decreasing sequence of positive integers $$\lambda=(\lambda_1\geq\lambda_2\geq \ldots\geq\lambda_l),$$
such that $\lambda_1+\lambda_2+ \cdots+\lambda_l=n.$ The $\lambda_i$'s are called the parts of $\lambda$ and $l$ is its size. The book \cite{A} is a classic in the theory
of integer partitions. A famous identity related to partitions and which plays an important role in  this paper is the  \textbf{First Rogers-Ramanujan Identity:}

\begin{center}

\textit{ The number of partitions of $n$ with no consecutive parts, neither equal parts is equal to the number of of partitions of $n$ whose parts are congruent to $1$ or $4$ modulo $5.$} 
 \end{center}
In \cite{BMS} (see also \cite{BMS1}) we got to this identity by considering the space of arcs  of the double point $X=\mbox{Spec} (\k[x]/x^2)$  centred at the origin; denote it 
by $X^0_\infty.$ The coordinate ring $A$ of $X^0_\infty$ is naturally graded and we  associate with it its Hilbert-Poincar\'e series 
that we call the Arc Hilbert-Poincar\'e series; we denote it by $AHP_{X,0}(t).$ Note that this is an invariant of singularities of algebraic varieties \cite{M,BMS}. We prove in \cite{BMS} that  
$AHP_{X,0}(t)$ is equal to the generating sequence of the number 
of partitions appearing in the Rogers-Ramanujan identities. The proof uses a Groebner basis computation associated with a monomial ordering (reverse lexicographical). Note 
that we have 
$$A=\frac{\k[x_i,i \in \textbf{Z}_{>0}]}{[x_1^2]},$$
where $[x_1^2]$ is the differential ideal generated by $x_1^2$ and its iterated derivative by the derivation
$D$ which is determined by $D(x_i)=x_{i+1}.$ So
$$[x_1^2]= (x_1^2,2x_1x_2,2x_1x_3+2x_2^2,\ldots).$$ 

The grading of $A$ is determined by giving to $x_i$ the weight $i.$ \\

In the frist part of this article, we consider a different monomial ordering (lexicographical) and we find that the Groebner basis computation with respect to this monomial
ordering is related with  an other identity involving the partitions that appear in the first Rogers-Ramanujan identity. We then use this other member in the Rogers-Ramanujan
identitiy to carry the computations of the Groebner-basis in small degrees (usual degree) but in all weights. This leads us to prove, in contrast with the case of the reverse
lexicographical ordering, that the Groebner basis of our ideal with respect to the lexicographical ordering is not differentially finite.\\

In the other part of this article, we will give a direct and simpler proof of a theorem by Nguyen Duc Tam \cite{N} where he computes a Groebner basis of the arc space of
$X=\mbox{Spec} (\k[x,y]/(xy))$ or simply of the differential ideal generated by $xy.$  We then use this theorem to obtain identities of partitions with two colors.

\section{Hilbert series and integer partitions}

In this section, we will begin by considering the Hilbert-Poincar\'e  series of some graded algebras that are inspired by Groebner basis computations; we then interpret these
series as generating sequences of partitions having special properties. At the end of the section, we will give more explanations on how we have found these graded 
algebras.\\ 

We will denote by $\k$ an algebraically closed field of characteristic $0.$ Recall that for a graded $\k-$algebras, $A=\oplus_{i \in \n}A_i$ (we assume that 
$\mbox{dim}_\k (A_i) <\infty$) the Hilbert-Poincar\'e series of $A,$ that
we denote by $HP(A),$ is by definition $$HP(A)=\sum_{i\in \n}\mbox{dim}_\k (A_i)q^i,$$
where $q$ is a variable.
For more about Hilbert-Poincar\'e series, see the appendix in \cite{BMS} and the 
references there.  

Let $n,l\geq 1$ be integer numbers. We consider the graded algebra $\k[x_l,x_{l+1},\ldots]$ where we give $x_i$ the weight $i,$ for 
$i\geq l.$  This grading induces a grading on the $\k$-algebra $\frac{\k[x_l,x_{l+1},\ldots]}{(x_{i_1}\cdots x_{i_n},~~i_j\geq l)}.$ \\

We consider the Hilbert-Poincar\'e series $$H_n^l=HP\left(\frac{\k[x_l,x_{l+1},\ldots]}{(x_{i_1}\cdots x_{i_n},~~i_j\geq l)}\right).$$

\begin{lem}\label{rec} We have
$$H_n^l=q^lH_{n-1}^l+H_n^{l+1}.$$

\end{lem}

\begin{proof} Using corollary 6.2 in \cite{BMS}, we have
$$H_n^l=HP\left(\frac{\k[x_l,x_{l+1},\ldots]}{(x_{i_1}\cdots x_{i_n},~~i_j\geq l)}\right)=$$
$$HP\left(\frac{\k[x_l,x_{l+1},\ldots]}{(x_l,x_{i_1}\cdots x_{i_n},~~i_j\geq l)}\right)+q^lHP\left(\frac{\k[x_l,x_{l+1},\ldots]}{\left((x_{i_1}\cdots x_{i_n},~~i_j\geq l): x_l\right)}\right)$$
The last term is exactly $H_n^{l+1}+q^lH_{n-1}^l.$

\end{proof}

\begin{prop}\label{hnl}
$$H_n^l=1+\frac{q^l}{1-q}+\frac{q^{2l}}{(1-q)(1-q^2)}+\cdots+ \frac{q^{(n-1)l}}{(1-q)(1-q^2)\cdots (1-q^{n-1})}.$$

\end{prop}
\begin{proof} The proof is by induction on the integer $n.$ Notice that for $n=1,$
$H_1^l=HP(\k)=1.$ For $n=2,$ the weighted-homogeneous components of
$$\frac{\k[x_l,x_{l+1},\ldots]}{(x_{i_1}x_{i_2},~~i_j\geq l)}$$
are generated by $1$ in degree $0$ and $x_i$ in degree $i$ for $i\geq l;$ 
for $i=1,\ldots l-1,$ the weighted-homogeneous components of degree  $i$ is the null vector space. Let us assume that the formula is true for $H_{j}^l,j\leq n-1$ and 
prove it for $H_n^l.$ Using lemma \ref{rec} repetitively, we obtain
$$H_n^l=q^lH_{n-1}^l+H_n^{l+1}=q^lH_{n-1}^l+q^{l+1}H_{n-1}^{l+1}+H_n^{l+2}=\cdots=$$
$$q^lH_{n-1}^l+q^{l+1}H_{n-1}^{l+1}+\cdots+q^mH_{n-1}^{m}+H_n^{m+1}.$$
But we have $$\lim_{m \to \infty} H_n^m = 1,$$
where the limit is considered for the $q$-adic topology in $\c[[q]];$ hence we can write
$$ H_n^l=1+ q^lH_{n-1}^l+q^{l+1}H_{n-1}^{l+1}+\cdots+q^mH_{n-1}^{m}+q^{m+1}H_{n-1}^{m+1}+\cdots=$$
$$1+\sum_{m\geq l}q^mH_{n-1}^{m}. $$
By the induction hypothesis we obtain
$$ H_n^l=1+$$
$$q^l+\frac{q^{2l}}{1-q}+\frac{q^{3l}}{(1-q)(1-q^2)}+\cdots+ \frac{q^{(n-1)l}}{(1-q)(1-q^2)\cdots (1-q^{n-2})}+$$
$$q^{l+1}+\frac{q^{2(l+1)}}{1-q}+\frac{q^{3(l+1)}}{(1-q)(1-q^2)}+\cdots+ \frac{q^{(n-1)(l+1)}}{(1-q)(1-q^2)\cdots (1-q^{n-2})}+
\cdots=$$
(summing by columns) 
$$1+(q^l+q^{l+1}+\cdots) + (\frac{q^{2l}}{1-q}+\frac{q^{2(l+1)}}{1-q}+\cdots)+\cdots +$$
$$(\frac{q^{(n-1)l}}{(1-q)(1-q^2)\cdots (1-q^{n-2})}+
\frac{q^{(n-1)(l+1)}}{(1-q)(1-q^2)\cdots (1-q^{n-2})}+\cdots).$$
The formula for $H_n^l$ follows from the formula
$$q^{jl}+q^{j(l+1)}+q^{j(l+2)}+\cdots=\frac{q^{jl}}{1-q^{j}},~j=1,\ldots,n-1.$$

\end{proof}

\begin{cor}Let $l,n\geq 1$ be integers. The generating series of the integer partitions with parts
greater or equal to $l$ and size (number of parts) less or equal to $n-1$ is
$$1+\frac{q^l}{1-q}+\frac{q^{2l}}{(1-q)(1-q^2)}+\cdots+ \frac{q^{(n-1)l}}{(1-q)(1-q^2)\cdots (1-q^{n-1})}.$$

\end{cor}
\begin{proof}
 This is just the combinatorial interpretation of proposition \ref{hnl}: A basis of the $i-$th weighted-homogeneuous component of 
 $\frac{\k[x_l,x_{l+1},\ldots]}{(x_{i_1}\cdots x_{i_n},~~i_j\geq l)}$ is given by the monomials $x_{h_1}\cdots x_{h_r}$ such that
 $h_1+\cdots+h_r=i,h_j\geq l$ and  $x_{h_1}\cdots x_{h_r} \not\in (x_{i_1}\cdots x_{i_n},~~i_j\geq l);$ this corresponds 
 to the partitions of $i$  with parts
greater or equal to $l$ and size (number of parts) less or equal to $n-1.$
\end{proof}

We now look at the Hilbert series that makes the link to the first Rogers-Ramanujan identity. For $j\geq 1,$ let 
\begin{equation}\label{H1}
 HP\left(\frac{\k[x_j,x_{j+1},\ldots]}{(x_{i_1}\cdots x_{i_k}x_k,~~i_1 \geq i_2\geq \ldots \geq i_k\geq k\geq j)}\right).
\end{equation}

\begin{lem}\label{recp} We have
$$H_j=q^jH^j_{j}+H_{j+1}.$$

\end{lem}

\begin{proof} Using corollary 6.2 in \cite{BMS}, we have
$$H_j=HP\left(\frac{\k[x_j,x_{j+1},\ldots]}{(x_{i_1}\cdots x_{i_k}x_k,~~i_1\geq i_2\geq\ldots\geq i_k\geq k\geq j)}\right)=$$
$$HP\left(\frac{\k[x_j,x_{j+1},\ldots]}{(x_j,x_{i_1}\cdots x_{i_k}x_k,~~i_1\geq i_2\geq\ldots\geq i_k\geq k\geq j)}\right)+$$
$$q^jHP\left(\frac{\k[x_j,x_{j+1},\ldots]}{((x_{i_1}\cdots x_{i_k}x_k,~~i_1\geq i_2\geq\ldots\geq i_k\geq k\geq j):x_j)}\right).$$
The last term is exactly $H_{j+1}+q^jH^j_{j}.$

\end{proof}

\begin{theo}\label{HP}We have
$$H_1=1+\frac{q}{1-q}+ \frac{q^4}{(1-q)(1-q^2)}+ \frac{q^9}{(1-q)(1-q^2)(1-q^3)}+\cdots=$$
$$1+\sum_{n\geq1}\frac{q^{n^2}}{(1-q)(1-q^2)\cdots(1-q^n)}.$$ 
\end{theo}
\begin{proof}Applying lemma \ref{recp} for $j=1$, we obtain that 
$H_1=qH_1^1+H_2;$ applying the same lemma for $j=2$ we have that $H_1=qH_1^1+q^2H^2_2
+H_3;$ Applying repetitively and in the same way  lemma \ref{recp}, we obtain that  
for $m\geq 2,$
$$H_1= qH^1_{1}+q^2H^2_{2}+\cdots+q^mH^m_{m}+H_{m+1}.$$
Noticing that $$\lim_{m \to \infty} H_m = 1$$ (where the limit is considered for the $q$-adic topology in $\c[[q]], \c$ being the field of complex numbers) we can write 
 $$H_1=1+qH^1_{1}+q^2H^2_{2}+\cdots+q^mH^m_{m}+\cdots.$$ Using proposition \ref{hnl}, we obtain that $H_1$ is equal to \\
 
 \xymatrix{
      1 & +q &  \\
       & +q^2  & +\frac{q^4}{1-q} & \\
      &  +q^3 & +  \frac{q^6}{1-q} & +\frac{q^9}{(1-q)(1-q^2)}\\
      &  +q^4 & +  \frac{q^8}{1-q} & +\frac{q^{12}}{(1-q)(1-q^2)} &+~~~\cdots\\ } 

 Summing by columns we obtain the result.
 
\end{proof}

An equivalent statement of the theorem the following.

\begin{theo}\label{RRNM1} Let $n\geq 1$ be a positive integer. 
The number of partitions of $n$ with size less than or equal to the smallest part is equal to the number  of partitions of $n$ without consecutive nor equal parts.
\end{theo}
\begin{proof}
 This follows from the known fact (\cite{A,BMS}) that the series obtained in theorem \ref{HP} is also the generating sequence of the  partitions 
 without consecutive nor equal parts.
\end{proof}

\begin{theo}\label{RRNM2} Let $n\geq k$ be a positive integer. 
The number of partitions of $n$ with parts larger or equal to $k$ and size less than or equal to (the smallest part minus $k-1$) is equal to the number of  partitions of $n$ with parts larger or equal to $k$ and without consecutive nor equal parts.
\end{theo}
\begin{proof}Let denote by $H_k$ the Hilbert series of 
$$\frac{\k[x_k,x_{k+1},\ldots]}{(x_{i_1}\cdots x_{i_{j-k+1}}x_j,~~i_1\geq i_2\geq\ldots\geq i_{j-k+1}\geq j\geq k)}.$$
By using the same method as used in the proof of the theorem \ref {HP} we have
$$H_k=1+q^k H^k_{1}+q^{k+1}H^{k+1}_{2}+\cdots+q^{k+m}H^{k+m}_{m+1}+\cdots.$$

By proposition \ref {hnl} we have that

 \xymatrix{
 H_k=     1 & +q^k &  \\
       & +q^{k+1}  & +\frac{q^{2k+2}}{1-q} & \\
      &  +q^{k+2} & +  \frac{q^{2k+4}}{1-q} & +\frac{q^{3k+6}}{(1-q)(1-q^2)}\\
      &  +q^{k+3} & +  \frac{q^{2k+6}}{1-q} & +\frac{q^{3k+9}}{(1-q)(1-q^2)} &+~~~\cdots\\ } 

summing by columns we obtain that

$$H_k=1+\sum_{n\geq 1} \frac{q^{n(n+k-1)}}{(1-q)(1-q^2)\cdots(1-q^n)}.$$
So $H_{k+1}+q^k H_{k+2}$ equals to

 \xymatrix{
    1 & +\frac{q^{k+1}}{1-q} & +\frac{q^{2k+4}}{(1-q)(1-q^2)}& +\frac{q^{3k+9}}{(1-q)(1-q^2)(1-q^3)} & +~~~\cdots\\
      &  +q^k & +  \frac{q^{2k+2}}{1-q} & +\frac{q^{3k+6}}{(1-q)(1-q^2)} &+~~~\cdots\\ } 

Again summing by columns we obtain that 

$$H_{k+1}+q^k H_{k+2}=1+\sum_{n\geq 1} \frac{q^{n(n+k-1)}}{(1-q)(1-q^2)\cdots(1-q^n)}=H_{k}.$$ 

Using this equation repetitively proves that for each $k\geq 1$  we can write following recursion series

$$H_k=A_{k+i}H_{k+i}+B_{k+i+1}H_{k+i+1}.$$

With initial conditions $A_{k+1}=1$ and $B_{k+2}=q^k$, and for all $i\geq 2$ we have $A_{k+i}=A_{k+i-1}+B_{k+i}$ and $B_{k+i+1}=q^{k+i-1}A_{k+i-1}.$

But if we denote by $H'_k$ the Hilbert series of $\frac{\k[x_k,x_{k+1},\ldots]}{(x_i^2,x_{i+1}^2,i\geq k)},$ then by \cite{BMS}, $H'_k$ satisfies in the following recursion formula

$$H'_k=A_{k+i}H'_{k+i}+B_{k+i+1}H'_{k+i+1}.$$
Where $A_j,B_j\in\k[[q]]$ are the same as before. So in the $(q)$-adic topology we have

$$H_k=\lim\big( A_{k+i}H_{k+i}+B_{k+i+1}H_{k+i+1}\big)=lim A_{k+i}=\lim\big( A_{k+i}H'_{k+i}+B_{k+i+1}H'_{k+i+1}\big)=H'_k.$$
This equation gives us the result. 
\end{proof}
Theorem \ref{RRNM1} is inspired from a Groebner basis computation of the differential ideal $[x_1^2]:$ By  \cite{BMS}, the initial ideal of $[x_1^2]$
with respect to the reverse lexicographical ordering is $(x_i^2,x_ix_{i+1},i\geq 1),$ while we can guess that its initial ideal 
with respect to the lexicographical ordering is $((x_{i_1}\cdots x_{i_k}x_k,~~i_1 \geq i_2\geq \ldots \geq i_k\geq k\geq 1)$ (we make use of this 
guess in the next section). Hence the Hilbert series of the quotient rings of  $\k[x_1,x_2,\ldots]$ by these ideals are equal. The Hilbert series
of the quotient by $(x_i^2,x_ix_{i+1},i\geq 1),$  is the generating series of the partitions  without consecutive nor equal parts; 
The Hilbert series
of the quotient by $((x_{i_1}\cdots x_{i_k}x_k,~~i_1 \geq i_2\geq \ldots \geq i_k\geq k\geq 1)$  is the generating series of the partitions  
 with size less or equal to the smallest part. Theorem \ref{RRNM2} can be guessed in the same way by considering the ideal 
 $[x_k^2]$ in $\k[x_k,x_{k+1},\ldots].$

\section{On the lex Groebner basis of $[x_1^2]$ }

Again, let $\k$ be a field of characteristic $0,$ and consider the 
graded ring $\k[x_1,x_2,\dots],$ where the weight of $x_i$ is  $i$. So the weight of the monomial $x^{\alpha}:= x_{i_1}^{\alpha_{1}} x_{i_2}^{\alpha_{2}} \dots x_{i_n}^{\alpha_{n}} $, is equal to $\sum_{j=1}^{n} i_j {\alpha_{j}} $, and its (usual) degree is equal to $n$.\\
Let $f_2=x_1^2$ and for $i\geq 3,f_i=D^{i-2}(f_2):=D(f_{i-1}),$ where $D$ is the derivation determined by $D (x_i)=x_{i+1};$  then $I=[f_2]:=(f_2,f_3,\dots)\subset  \k[x_1,x_2,\dots],$ is the defining ideal  (up to isomorphism) of the space of arcs centred at the origin of  $X=\mbox{Spec}(\k[x]/(x^2)).$ The ideal $I$ is a differential ideal, i.e. we have $D(I)\subset I.$ We are interested in the possibility that $I$ have a differentially finite Groebner basis with respect to a monomial ordering ; see the following definition. \\ 

\begin{defi} Let $J\subset \k[x_1,x_2,\dots]$ be a differential ideal with respect to $D$ (i.e. $D(J)\subset J).$ Let $"<"$ be a total monomial order defined on $\k[x_1,x_2,\dots].$
We say that $J$ has a differentially finite Groebner basis with respect to $"<",$ if  there exist a finite number of polynomials $h_1,\ldots,h_r \in \k[x_1,x_2,\dots]$ such that $J=[h_1,\ldots,h_r]$ and the initial ideal $In_<(J)$ of $J$ with respect to $"<"$ satisfies
$$ In_<(J)=(In_<(D^i(h_j)),j=1,\ldots,r;i\geq 0),$$
where $D^i$ denotes the $i-$th iterated derivative and $D^0$ is the identity.

\end{defi}

Note that the notation $[h_1,\ldots,h_r]$ in the definition denotes the differential ideal generated by the $h_i, i=1,\ldots,r$ and by all their iterated derivatives. . Note that there might be different notions of differential Groebner basis, see \cite{CF}, \cite{O} and their bibliography.\\

In this section, we prove that no Groebner basis of $I$ with respect to 
the weighted lexicographical order is  differentially finite. Note that, in contrary
with this case, it follows for \cite{BMS} that in the case of the weighted reverse lexicographical order $I$ has a differentially finite Groebner basis. \\

We denote the $n$-th derivative of a polynomial $f_i$ by $f_i^{(n)}$, so we have

$f_n=f_2^{(n-2)}=(x_1^2)^{(n-2)}=\sum\limits_{i=0}^{n-2} \binom {n-2}{i} x_1^{(i)}x_1^{(n-2-i)}=\sum\limits_{i=0}^{n-2} \binom {n-2}{i} x_{1+i}x_{n-i-1}.$
 
 Denote the leading term of $f_n$ with respect to the weighted lexicographical order by $LT(f_n).$ So $LT(f_n)=2 x_1 x_{n-1}$ for all $n\geq 2.$ \\
 
Recall that the S-polynomial of $f,g \in  \k[x_1,x_2,\dots]$ is by definition 
$$S(f,g):=\frac{x^\gamma}{LT(f)}f-\frac{x^\gamma}{LT(g)}g,$$
where $x^\gamma$ is the least common multiple of the leading monomials of $f$ and $g.$
A possible reference about S-polynomials and Groebner basis is \cite{GP}.\\

A direct computation of  the S-polynomial of $f_3$ and $f_4$ gives
$$S(f_3,f_4)=x_2^3.$$ We set  $F_{x_2^3}:=S(f_3,f_4).$ For $k>2,$
we recursively define 
$$F_{x_2x_k^2}:=S(F_{x_2 x_{k-1} ^2}^{(2)} , S(f_{k+1},f_{k+2})).$$
We then have the following lemma.

\begin{lem}\label{lead}
With respect to the weighted lexicographic order, for $k>2,$ 
the leading monomial of $F_{x_2x_k^2}$ is $x_2x_k^2.$
\end{lem}

\begin{proof}The proof is by induction on the integer $k$. Notice that for $k=3$ we have $F_{x_2^3}^{(1)}=3x_2^2  x_3$, $F_{x_2^3}^{(2)}=3x_2^2  x_4 + 6x_2  x_3^2$ and $S(f_4,f_5)=x_2^2x_4 - 3x_2x_3^2$. So we have
 $$ S(F_{x_2^3}^{(2)},S(f_4,f_5))=5x_2 x_3^2  :=F_{x_2 x_3^2}.$$
 For  $k=4$ we have $F_{x_2 x_3^2}^{(1)}=10 x_2 x_3 x_4+5 x_3^3$, $F_{x_2 x_3^2}^{(2)}=10 x_2 x_3 x_5+10 x_2 x_4^2 +25 x_3^2 x_4$ and $S(f_5,f_6)=3 x_2 x_3 x_5- 4x_2 x_4^2 -3 x_3^2 x_4$. So
 $$S(F_{x_2 x_3^2}^{(2)},S(f_5,f_6))=\frac{7}{3} x_2 x_4^2 +\frac{7}{2} x_3^2 x_4:=F_{x_2 x_4^2}.$$
Now assume that claim holds for $k-1\geq 4$. This means that for $k-1\geq 4$ we have 

$$F_{x_2 x_{k-1} ^2}:=S(F_{x_2 x_{k-2} ^2}^{(2)} , S(f_{k},f_{k+1}) ).$$
Since the leading monomial of $F_{x_2 x_{k-1} ^2}$ is $x_2 x_{k-1} ^2$, we can assume that  $F_{x_2 x_{k-1} ^2}=a x_2 x_{k-1}^2 + g_3$, for some rational number $a$ and some  polynomial $g_3$ with the monomials of the form $ x_{i_1} x_{i_2} x_{i_3}$ such that $3\leq i_1 \leq i_2 \leq i_3$. So, on one hand, the second derivative of $F_{x_2 x_{k-1} ^2}$ will be as follow
$$F_{x_2 x_{k-1} ^2}^{(2)}=2a x_2x_{k-1}x_{k+1}+2a x_2 x_k^2 +h_3.$$
Where $h_3=4a x_3x_{k-1}x_k +a x_4x_{k-1}^2 +g_3^{(2)}$. On the other hand, we have

$S(f_{k+1},f_{k+2})=S(\sum\limits_{i=0}^{k-1} \binom {k-1}{i} x_{1+i}x_{k-i},\sum\limits_{i=0}^{k} \binom {k}{i} x_{1+i}x_{k+1-i})$

$=\frac{1}{2}\sum\limits_{i=0}^{k-1} \binom {k-1}{i} x_{1+i}x_{k-i} x_{k+1}-\frac{1}{2}\sum\limits_{i=0}^{k} \binom {k}{i} x_{1+i}x_{k+1-i} x_k$

$=\frac{1}{2}\sum\limits_{i=1}^{k-2} \binom {k-1}{i} x_{1+i}x_{k-i} x_{k+1}-\frac{1}{2}\sum\limits_{i=1}^{k-1} \binom {k}{i} x_{1+i}x_{k+1-i} x_k$.

So by the above equation we obtain that $LT(S(f_{k+1},f_{k+2}))=(k-1)x_2 x_{k-1} x_{k+1}$.
\\Now we can compute $S(F_{x_2 x_{k-1} ^2}^{(2)} , S(f_{k+1},f_{k+2}))$.

$S(F_{x_2 x_{k-1} ^2}^{(2)} , S(f_{k+1},f_{k+2}) )$
\\$=\frac{1}{2a}(2a x_2 x_{k-1} x_{k+1}+ 2a x_2 x_k^2 + h_3)- \frac{1}{(k-1)}(\frac{1}{2}\sum\limits_{i=1}^{k-2} \binom {k-1}{i} x_{1+i}x_{k-i} x_{k+1}-\frac{1}{2}\sum\limits_{i=1}^{k-1} \binom {k}{i} x_{1+i}x_{k+1-i} x_k)$
\\$=x_2 x_k^2 +\frac{1}{2a} h_3- \frac{1}{(k-1)}(\frac{1}{2}\sum\limits_{i=2}^{k-3} \binom {k-1}{i} x_{1+i}x_{k-i} x_{k+1}-\frac{1}{2}\sum\limits_{i=1}^{k-1} \binom {k}{i} x_{1+i}x_{k+1-i} x_k)$
\\$=x_2 x_k^2 +\frac{1}{2a} h_3+\frac{k}{k-1} x_2 x_k^2- \frac{1}{(k-1)}(\frac{1}{2}\sum\limits_{i=2}^{k-3} \binom {k-1}{i} x_{1+i}x_{k-i} x_{k+1}-\frac{1}{2}\sum\limits_{i=2}^{k-2} \binom {k}{i} x_{1+i}x_{k+1-i} x_k)$
$=\frac{2k-1}{k-1}x_2 x_k^2+\frac{1}{2a} h_3- \frac{1}{(k-1)}(\frac{1}{2}\sum\limits_{i=2}^{k-3} \binom {k-1}{i} x_{1+i}x_{k-i} x_{k+1}-\frac{1}{2}\sum\limits_{i=2}^{k-2} \binom {k}{i} x_{1+i}x_{k+1-i} x_k).$

In the first sum, $2\leq i \leq k-3$ and in the second one $2\leq i\leq k-2$. So each monomial that appears in $S(F_{x_2 x_{k-1} ^2}^{(2)} , S(f_{k+1},f_{k+2}))$ is of the form $x_3 x_{i_1}  x_{i_2}  x_{i_3}$ such that $3\leq i_1 \leq i_2 \leq i_3$, except $x_2 x_k^2$ and hence 
$$LT(S(F_{x_2 x_{k-1} ^2}^{(2)}, S(f_{k+1},f_{k+2})))=\frac{2k-1}{k-1}x_2x_k^2.$$ 
\end{proof}

\begin{theo}\label{NDF}  A Groebner  basis of the ideal $I,$ with respect to the weighted lexicographic order, is not  differentially finite.
\end{theo}
\begin{proof}

 For proving this fact, we will use the idea of Buchberger's algorithm (mainly that any cancellation of initial monomials comes from an $S-$ polynomial \cite{CLO}) to construct a part of a Groebner  basis of the ideal $I$ with respect to the weighted lexicographic order, which is differentially infinite.\\
 
By lemma \ref{lead}, we have that for every integer $n\geq 3,$ the initial monomial of the polynomial $F_{x_2x_n^2}$  is included in the initial ideal of $I.$\\ 

Let $G= \{f_i,F_{x_2x_n^2}, F_{x_2x_n^2}^{(m)} | i \geq 2, m\geq 1, n\geq 3 \} $. By Buchberger's algorithm $G$ may be a part of a Groebner  basis of the ideal $I$ but it is not a Groebner  basis of $I$  because 

$$S(F_{x_2x_3^2},F_{x_2x_3^2}^{(1)})=S(F_{x_2x_3^2},F_{x_2x_3x_4})
=5x_3^4.$$

But the monomial $x_3^4$ which is a member of the ideal $I$ is not divisible by the leading terms of any element of $G.$ \\

Note that the (usual) degree of the $S-$polynomial of two polynomials is at least equal to the maximum of degrees of these two polynomials. On the other hand, the derivative of a polynomial has the same degree as itself. Hence 
the  degree of the $f_i$'s is equal to two, and other elements of $G$ have degree strictly bigger than $2.$
\\This means that the monomials of degree two that appear as the leading terms of elements of Groebner basis, are of the form $x_1x_i$ for $i\geq1$.
\\So we do not have any polynomial in a Groebner basis whose leading monomial is $x_2 x_n $ for some $n\geq 2$, and so  a polynomial having the same initial monomial of $F_{x_2x_n^2}$ should be included in this Groebner 
basis for each integer $n\geq 3.$ Since the initial monomial of the polynomial $F_{x_2x_n^2}$ is not the initial of the derivative of any other element of $G,$ the Groebner  basis of the ideal $I$ with respect to the weighted lexicographic order will not be 
differentially finite: it should contain  polynomials whose initial monomials are the initials of $F_{x_2x_n^2}, n \geq 3$ and no one of these initial monomials is the derivative of an other initial monomial of an element in $G.$

\end{proof}


\section{Two colors partitions and the node}

Let $S:=\k[x_1,x_2,\ldots,y_1,y_2,\dots]$ be the graded polynomial ring where
$x_i,y_i$ have the weight $i$ for every $i\geq1;$ the order of appearance of the variables
is important since we will use below a reverse lexicographical ordering. We consider the derivation
on $S$ defined by $D(x_i)=x_{i+1}$ and $D(y_i)=y_{i+1}.$ Let $f_2=x_1y_1,$ and let
$$I=[f_2]=(x_1y_1,x_2y_1+x_1y_2,\ldots)$$
be the ideal generated by $xy$ and its iterated derivatives $f_i,i\geq 3$ by $D:$
for $i\geq 3, f_i=D(f_{i-1}).$  Note that the scheme defined by $I$ is the space of arcs centred at the origin of the node $X=\{xy=0\} \subset \textbf{A}^2.$\\ 
In this section, we are interested in determining a Groebner basis of $I$ with
respect to the weighted reverse lexicographical order and then to apply this result
to integer partitions. Note that the Groebner basis below was found by Nguyen Duc Tam \cite{N}; he has a beautiful but very long and difficult proof that this is  actually a Groebner basis. Below we give a simpler 
and very short proof. 
\\We will begin by defining elements of $I,$ and we will show later that these elements give the Groebner basis cited above.


\begin{defi}(\cite{N})
For $1 \leq i_1 \leq i_2 \leq \dots \leq i_k $ and for $k\geq 2,$ we set 
$$G_{i_1 ,i_{2}+1,i_{3}+2,\dots ,i_{k}+k-1}:=
\det{\begin{bmatrix}
 x_{i_1 -k+2}&x_{i_1 -k+3}&\dots &x_{i_{1}}& f_{i_1 +1}\\
x_{i_2 -k+3}&x_{i_2 -k+4}&\dots &x_{i_{2}+1}& f_{i_2 +2}\\
x_{i_3 -k+4}&x_{i_3 -k+5}&\dots &x_{i_{3}+2}& f_{i_3 +3}\\
\vdots & \vdots & \ddots & \vdots &\vdots\\
x_{i_k +1}&x_{i_k +2}&\dots &x_{i_{k}+k-1}& f_{i_k +k}
\end{bmatrix}}$$
where det stands for determinant.  
\end{defi}
Expanding the determinant with respect to the last column, we see that these are elements of $I.$ A direct computation using the definition of the $f_i$ gives the 
following:
\begin{lem}\label{IN}\cite{N}
The leading term of  $G_{i_1 ,i_{2}+1,i_{3}+2,\dots ,i_{k}+k-1}$ with respect to weighted reverse lexicographic order, is $x_{i_1} x_{i_2} x_{i_3}\dots x_{i_k}y_k$.
\end{lem}
 We denote by $\mathbb{G}$ the set whose elements are the $ G_{i_1 ,i_{2}+1,i_{3}+2,\dots ,i_{k}+k-1}$ and the $f_i.$ It follows from lemma \ref{IN} that the ideal generated
 by the initials of the elements of $\mathbb{G}$ is
$$J:=(x_{i_1} x_{i_2}\dots x_{i_k} y_k |~~ i_j ,k \geq 1). $$
 


First, we are interested in computing the Hilbert-Poincar\'e series of $S/J.$ For that we introduce for $n\geq 1$  the Hilbert-Poincar\'e series 
$$HP_n=HP\left(\frac{\k [x_i,y_j | i\geq 1 , j\geq n]}{(x_{i_1} x_{i_2}\dots x_{i_k} y_k |~~ i_j\geq 1 ,k \geq n)}\right).$$ So $HP(S/J)=HP_1.$ We will use the following form of $H_n^1$ from
section 1:

\begin{lem}\label{HN1} For any  $n\geq 2$ we have
$$H^1_n = \frac{1}{(1-q)(1-q^2)\dots (1-q^{n-1})}.$$
\end{lem}

\begin{proof}
By proposition $\ref{hnl}$ we have
$$H^1_n = 1+ \frac{q}{1-q} +\frac{q^2}{(1-q)(1-q^2)} +\dots + \frac{q^{n-1}}{(1-q)\dots(1-q^{n-1})}$$
We prove the expression in the lemma by induction on the integer $n$. For $n=2$
$$H^1_2 = 1+\frac{q}{1-q} =\frac{1}{1-q}$$
Assume that $H^1_n= \frac{1}{(1-q)\dots (1-q^{n-1})}$.
Now  we have
$$H^1_{n+1} = 1 +\frac{q}{1-q} + \dots +\frac{q^{n-1}}{(1-q)\dots (1-q^{n-1})}+ \frac{q^n}{(1-q)\dots (1-q^{n})}$$
$$=H^1_n +\frac{q^n}{(1-q)\dots (1-q^n)}$$
By induction hypothesis we obtain
$$H^1_{n+1} = \frac{1}{(1-q)\dots (1-q^{n-1})} +\frac{q^n}{(1-q)\dots (1-q^n)}=\frac{1}{(1-q)\dots (1-q^n)}$$
\end{proof}

\begin{lem}\label{rec1}
$$HP_n =HP_{n+1} + q^n \prod_{i\geq1} \frac{1}{1-q^i}.$$
\end{lem}
\begin{proof}Using corollary 6.2 in \cite{BMS}, we have
$$HP_n=HP\left (\frac{\k [x_i,y_j | i\geq 1 , j\geq n]}{(x_{i_1} x_{i_2}\dots x_{i_k} y_k |~~ i_j \geq 1,k \geq n)}\right)=$$
$$HP\left(\frac{\k [x_i,y_j | i\geq 1 , j\geq n+1]}{(x_{i_1} x_{i_2}\dots x_{i_k} y_k |~~ i_j\geq 1 ,k \geq n)}\right)+q^n HP\left(\frac{\k [x_i,y_j | i\geq 1 , j\geq n]}{(x_{i_1} x_{i_2}\dots x_{i_n}  | ~~i_j \geq 1)}\right)=$$
$$HP_{n+1} + q^n HP\left(\frac{\k [x_1,x_2,\cdots]}{(x_{i_1} x_{i_2}\dots x_{i_n}  | ~~i_j \geq 1)}\right) ~~HP\left(\k[y_{n},y_{n+1},\cdots]\right)= $$
$$HP_{n+1} +q^n H_n^1 \prod_{i\geq n} \frac{1}{1-q^i} ;$$
by lemma \ref{HN1} this is equal to
$$HP_{n+1} +\frac{q^n}{(1-q)\cdots (1-q^{n-1})} \prod_{i\geq n} \frac{1}{1-q^i}= $$
$$HP_{n+1} +q^n \prod_{i\geq 1} \frac{1}{1-q^i}.$$
\end{proof}
\begin{prop}\label{HPn} We have
$$HP(S/J)=HP_{1}=\frac{1}{1-q}~~ \prod_{i\geq 1} \frac{1}{1-q^i}.$$

\end{prop}

\begin{proof}

Using lemma \ref{rec1} repetitively we obtain that, for $m\geq 2$
$$HP_1= q  \prod_{i\geq 1} \frac{1}{1-q^i} + q^2  \prod_{i\geq 1} \frac{1}{1-q^i} +\cdots q^m  \prod_{i\geq 1} \frac{1}{1-q^i}+ HP_{m+1}. $$
On the other hand
$$\lim_{m \to \infty} HP_m =  \prod_{i\geq 1} \frac{1}{1-q^i},$$
where the limit is considered for the $q$-adic topology in $\c[[q]];$ so we have,
$$HP_1 =  \prod_{i\geq 1} \frac{1}{1-q^i} + q  \prod_{i\geq 1} \frac{1}{1-q^i} + q^2  \prod_{i\geq 1} \frac{1}{1-q^i} + \cdots$$
$$=(1+q+q^2+\cdots )\prod_{i\geq 1} \frac{1}{1-q^i}$$
$$=\frac{1}{1-q}~~ \prod_{i\geq 1} \frac{1}{1-q^i}.$$

\end{proof}
We now are ready to prove:

\begin{theo}(\cite{N}) We have that $\mathbb{G}$ is a Groebner basis of $I.$

\end{theo}

    \begin{proof}
    Let $In(I)$ be the initial ideal of $I$ with respect to the weighted reverse
lexicographical order. Since all the elements of $\mathbb{G}$ are also in $I,$
we have that $J\subset In(I);$ to prove that $\mathbb{G}$ is a Groebner basis of $I,$ we need to prove that $J=In(I).$\\ 
Noticing that $(f_2,f_3,\ldots)$ is a regular sequence \cite{GS} (Note that this is rarely the case \cite{M1}) and that $f_i$ is of weight $i,$ we deduce that 
$$HP(S/I)=\frac{1}{1-q}~~ \prod_{i\geq 1} \frac{1}{1-q^i},$$
which is equal by proposition  \ref{HPn} to $HP(S/J).$ But since we have a flat
deformation with generic fibre $S/I$ and special fibre $S/In(I),$ we have
that $HP(S/I)=HP(S/In(I)),$ hence $HP(S/In(I))=HP(S/J).$ We deduce 
that the homogeneous components of the same weight of $S/(In(I)$ and $S/J$ have the same (finite) dimension, and since we have an inclusion in one sense because $J\subset In(I),$ they
are equal. Hence $J=In(I).$

  \end{proof}

We will interpret the above results in terms of two color partions: consider that we have two copies of each positive integer
number $m,$ one is blue and the other is red; we denote these copies by $m_b$ and $m_r.$   We define an order between the colored integers by $m_b>m_r$ (so that we do
not count in a partition  $m_b+m_r$ and $m_r+m_b$ as different); if $m>k,$ we say $m_c>k_{c'}$ for $c,c'\in\{b,r\}.$\\

An integer partition of a positive integer number $n$ is a decreasing sequence (with respect to the order that we have just defined) of positive integers of one color or an other
$$\lambda=(\lambda_{1,c_1}\geq\lambda_{2,c_2}\geq \ldots\geq\lambda_{l,c_l}),$$
where $c_i \in\{b,r\}$ and such that $\lambda_{1,c_1}+\lambda_{2,c_2}+ \cdots+\lambda_{l,c_l}=n.$
For example, the two colors integer partitions of $2$ are:
\begin{center}
 $2_b$\\
 $2_r$\\
 $1_b+1_b$\\
 $1_r+1_r$\\
 $1_b+1_r.$

\end{center}

Colored partitions has already appeared in the work of Andrews and Agarwal \cite{APS}.\\

On one hand, we can interpret the series $$\frac{1}{1-q}~~ \prod_{i\geq 1} \frac{1}{1-q^i}$$
as the generating series of the partitions with $2$ colors of $1$ and only the red color of any other positive integer. So the patitions of $2$ of this type
are all the partitions appearing in the above example except the first one.\\

On the other hand, the monomials in $S/J$ of weight $n$ are in bijection with the partitions with $2$ colors of $n$ whose number of blue parts is strictly less 
than its smallest red part (if this latter exists). In the above example of partions of $2,$ all the partitions but the last one are of this type.
The Hilbert-Poincar\'e series $HP(S/J)$ is then the generating sequence of this type of partitions. Hence proposition \ref{HPn} gives:

\begin{theo}\label{2c}
The number of partitions of $n$  with $2$ colors of $1$ and only the red color of any other positive integer is equal to the number of
partitions with $2$ colors of $n$ whose number of blue parts is strictly less 
than its smallest red part (if this latter exists).
\end{theo}

Playing the same game with the ideal $[x_jy_j]$ instead of $[x_1y_1]$ we can prove the following generalization of theorem \ref{2c}:

\begin{theo}
Let $j$ be a positive integer number. The number of partitions of $n$  with $2$ colors of $1,\ldots,2j-1$ and only the red color of any other positive integer is equal to the number of
partitions with $2$ colors of $n$ whose number of blue parts is strictly less 
than its smallest red part (if this latter exists) minus $(j-1).$
 
\end{theo}
We recover theorem \ref{2c} by putting $j=1.$

\newpage

Pooneh Afsharijoo$^1$, Hussein Mourtada$^2$\\

\noindent Equipe G\'eom\'etrie et Dynamique, \\
Institut Math\'ematique de Jussieu-Paris Rive Gauche,\\
 Universit\'e Paris Diderot, \\
 B\^atiment Sophie Germain, case 7012,\\
75205 Paris Cedex 13, France.\\

\noindent Email 1: pooneh.afsharijoo@imj-prg.fr\\
   Email 2: hussein.mourtada@imj-prg.fr

\end{document}